\documentclass[english]{article}
\usepackage[T1]{fontenc}
\usepackage[utf8]{inputenc}
\usepackage{textcomp}
\usepackage{amsmath}
\usepackage{amssymb}

\makeatletter
\newenvironment{lyxlist}[1]
{\begin{list}{}
{\settowidth{\labelwidth}{#1}
 \setlength{\leftmargin}{\labelwidth}
 \addtolength{\leftmargin}{\labelsep}
 }}
{\end{list}}

\usepackage{ifxetex}\ifxetex
\usepackage{fontspec}
\else
\fi

\usepackage{fixltx2e}\usepackage{fancyhdr}

\makeatother

\usepackage{babel}
\begin{document}

\title{\textbf{The Gödel's legacy: revisiting the Logic}}

\author{Dr. Giuseppe Raguní UCAM University of Murcia - Spain - graguni@ucam.edu }
\maketitle
\begin{abstract}
Some common fallacies about fundamental themes of Logic are exposed:
the First and Second incompleteness Theorem interpretations, Chaitin's
various superficialities and the usual classification of the axiomatic
Theories in function of its language order.\end{abstract}
\begin{lyxlist}{00.00.0000}
\item [{KEYWORDS:}] Incompleteness, undecidability, semantic completeness,
categoricity, randomness, Chaitin's constant, first and second order
languages, consistency.\end{lyxlist}
\begin{enumerate}
\item Prologue
\item The pressure of a sentence
\item Chaitin's licenses
\item A tired classification
\item About the interpretation of the Second incompleteness Theorem
\end{enumerate}

\section*{1. Prologue}

Regularly, after having enjoyed the fruits of the genius of an extraordinary
man, we have to suffer the dogmatism of any his conclusions: everything
he said is gold. In physics, an experiment can eliminate the most
stubborn opinion, but if the matter is in relation with pure philosophy,
things are much more complicated and rarely solvable in a compulsory
manner.

However, in some sectors \texttt{\textquotedbl{}}applied\texttt{\textquotedbl{}}
of philosophy, such as in mathematical Logic, since some time is possible
- and imperative - to require accuracy and rigor. Hilbert's formalism
introduced by the late nineteenth century, in fact, has equipped the
axiomatic Disciplines with a orderly symbolism, deprived of the ambiguity
of any intuitive intervention and regimented to a rigorous epistemological
analysis. As a consequence, the modern Logic has been capable of remarkable
results, led by the astonishing Gödel's theorems.

Despite the above, still today several serious errors are common in
this matter. These lapses, combined with the use of a dated terminology
which today is indubitably insidious and insufficient, are hindering
the spread of this invaluable and fascinating knowledge, and not only
in the humanistic field.

Here, I will try to clarify briefly some of these confusions, citing
some text where it is possible to find more depth.

\section*{2. The pressure of a sentence}

The first point concerns the range of applicability of the Gödel's
First incompleteness Theorem. We omit, now, both the statement and
the very important consequences of this famous theorem, that the reader
can find in any good book on Logic. The thing to emphasize here is
just that the same hypothesis of the Theorem require the \textit{enumerability}
of the set of theorems and proofs of the Theory to which it can be
applied%
\footnote{ A set is said \textit{enumerable} if it exists a biunivoc correspondence
between its elements and the natural numbers.%
}. As a preliminary to the demonstration, Gödel established a special
numerical code (called briefly \textit{gödelian}) both for propositions
and demonstrations. When applied to the \emph{formal}%
\footnote{ Throughout the article, with the adjective \texttt{\textquotedbl{}}formal\texttt{\textquotedbl{}}
we wish to mean \texttt{\textquotedbl{}}empty of explicit meaning\texttt{\textquotedbl{}},
or \texttt{\textquotedbl{}}symbolic\texttt{\textquotedbl{}}, \texttt{\textquotedbl{}}syntactic\texttt{\textquotedbl{}},
\texttt{\textquotedbl{}}codified\texttt{\textquotedbl{}}.%
} Theory of natural numbers (or Peano Arithmetic, \textit{PA} from
now), this encoding makes every proposition and every demonstrations
be assigned a numeric code, unique and exclusive. But what's happen
if you do the same with an arithmetic Theory that has got a number
\textit{not enumerable} of theorems? This Theory exists and is usually
called \texttt{\textquotedbl{}}second order Arithmetic\texttt{\textquotedbl{}}.
In its premises there is a meta-mathematic axiomatic scheme that,
generalizing the principle of induction of \textit{PA}, introduces
one axiom for each element of \textit{P(N)}, the set of all the subsets
of natural numbers. We will call briefly \textit{complete} this induction.
Since this set is, as well known, \textit{not enumerable}, it follows
that also the sets of theorems and proofs are not enumerable. Therefore,
in this Theory, not all the proofs can have got a different \textit{gödelian}:
else its would be enumerable. The not-denumerability of the proofs
reveals the indispensable use of a \textit{intrinsically semantic}
strategy (ie with use of not fully codifiable meanings) to derive
the theorems.

Even if we want to consider the principle of \textit{complete} induction
as one symbolic formula, ie as a single semantic axiom, this axiom
is \textit{not} decidable (or just effectively enumerable), since
its semanticity is not removable. This last conclusion is reached,
for example, representing the Theory inside the formal Set Theory%
\footnote{ Anyone between NBG, ZF or MK.%
}: the axiom has to be translated to an axiomatic scheme that generates
a number \textit{not enumerable} of inductive axioms%
\footnote{ Indeed, a common mistake is to believe that a finite number of statements
is always \textit{decidable}. More in: G. Raguní,\textit{ Confines
lógicos de la Matemática}, revista cultural \textit{La Torre del Virrey}
\textit{- Nexofía}, free on the \textit{Web}: http://www.latorredelvirrey.es/nexofia/pdf/Confines-logicos-de-la-matematica.pdf
(2011), pp. 153-158 and 293-295.%
}.

The simple consequence is that the First incompleteness Theorem can
not be applied to the second order Arithmetic%
\footnote{ In fact in this Theory, the statement of Gödel, whose standard interpretation
is \texttt{\textquotedbl{}}there is no code for a proof of this statement\texttt{\textquotedbl{}}
is not equivalent to \texttt{\textquotedbl{}}I am not provable in
this Theory\texttt{\textquotedbl{}}, as in \textit{PA}.%
}. However, this ambiguity has been long and loudly spread everywhere,
even in specifically technical areas%
\footnote{ Just two examples: \texttt{<\textcompwordmark{}<}\textit{ovviamente
il} \textit{primo teorema d'incompletezza è dimostrabile anche nell'Aritmetica
al secondo ordine}\texttt{>\textcompwordmark{}>}, E. Moriconi, \textit{I
teoremi di Gödel}, SWIF (2006), on the \textit{Web}: http://lgxserve.ciseca.uniba.it/lei/biblioteca/lr/public
/moriconi-1.0.pdf, p. 32; C. Wright, \textit{On Quantifying into Predicate
Position: Steps towards a New(tralist) Perspective} (2007), p. 22.
In this last work, maybe it is significant that the author discusses
this property with a number of delicate epistemologic questions. Anyway,
in both cases the property is considered \texttt{\textquotedbl{}}obvious\texttt{\textquotedbl{}}
without further explanation. %
}. About this subjet, often is perceived a typical strange lack of
rigor (perhaps ill-concealed hint of doubt). For example is repeated
usually that \texttt{\textquotedbl{}}also the second order Arithmetic,
since it contains the axioms of \textit{PA,} is subdued by the incompleteness
Theorem\texttt{\textquotedbl{}}. Forgetting that also must be kept
the \textit{effective axiomatizability.} Emblematic is the case of
the Theory obtained from \textit{PA} by adding, as axioms, all the
statements true in the \textit{standard model} of \textit{PA}: it
also contains all the axioms of \textit{PA,} but it is complete.

Really, it is even possible that the second order Arithmetic is syntactically
complete, although its \emph{language} is semantically incomplete.
Since this Theory is not-formal, the answer to this interesting question
could come only from the Metamathematics%
\footnote{ G. Raguní, \textit{op.} \textit{cit.} (footnote n. 4), p. 296-297
.%
}.

There is no doubt that this situation is due to the overreverence
for the Gödel's figure. In the presentation of his incompleteness
Theorem in the Congress of Königsberg (1930), Gödel announced the
result as \texttt{\textquotedbl{}}a proof of the semantic incompleteness
of Arithmetic, since it is categorical\texttt{\textquotedbl{}}%
\footnote{ K. Gödel, \textit{Collected Works I: publications 1929-1936}, eds.
S. Feferman et al., Oxford University Press (1986), p. 26-29.%
}. He speaks, undoubtedly, of the second order Arithmetic which, unlike
\textit{PA}, is categorical%
\footnote{ A Theory is said \textit{categorical} if it admits a single model
up to isomorphism. In simpler words (but also more inaccurate), if
it has got a single correct interpretation.%
}. Gödel, therefore, based his claim by applying - wrongly - the First
Incompleteness Theorem to the second order Arithmetic: in fact, semantic
incompleteness derives from syntactic incompleteness and categoricity.
Despite the error, the conclusion is correct: as a result of Löwenheim-Skolem
Theorem%
\footnote{ We include here, both the \texttt{\textquotedbl{}}up\texttt{\textquotedbl{}}
and \texttt{\textquotedbl{}}down\texttt{\textquotedbl{}} version.%
}, categoricity and infinity of a model are sufficient to ensure semantic
incompleteness of the language of any arbitrary Theory. However, before
1936, when it became popular the generalization of Malcev, no logician,
including Skolem and Gödel, realized all the puzzling consequences
of this important Theorem.

Now, being the greatest logician of all time, the question is not
only what led him to error, but also why his claim has never been
subsequently corrected. It is not easy to answer to the first question.
It should be noted that Gödel, at least up to 1930, rarely distinguishes
the two types of arithmetic Theories, so logically different each
from one. Indeed in that period, neither he nor any other logician
sometimes highlights the \textit{intrinsic} \textit{semanticity} of
Theories with an uncountable number of statements. And the consequences
of it.

Everyone can well understand, on the other hand, why we don't have
got any kind of correction about the above announcement. In fact it
involves an argument - semantic completeness / incompleteness and
its relation with categoricity - that basically, as it was afterward
understood, has nothing to do with the syntactic completeness and
therefore with the First incompleteness Theorem. During the 30's,
this subject was very topical, especially because of Hilbert's concerns
about the categoricity. Gödel began showing that the formal classical
language%
\footnote{ We prefer to use this expression rather than the more usual \texttt{\textquotedbl{}}First
Order Classical Logic\texttt{\textquotedbl{}} for reasons that will
become clear in paragraph 4.%
} is semantically complete: if it is consistent, always has got at
least one model%
\footnote{ Semantic Completeness Theorem, 1929.%
}. It followed that any formal classical Theory syntactically incomplete,
ie, with at least one undecidable statement \textit{I}, could not
be categorical. Because it supports at least two nonisomorphic models:
one for which \textit{I} is true and another for which is false. In
addition, Gödel continued believing - as tacitly assumed by Hilbert
and any other logician at time - that the syntactic completeness of
a formal Theory (or, more generally, of a Theory with a semantically
complete language) should imply its categoricity. Consequently the
First incompleteness Theorem was seen, just after its acceptance,
as an instrument capable of discriminating the categorical nature
and / or the semantic completeness of the axiomatic Systems. For example,
the formal arithmetic Theory, ie \textit{PA}, was believed to be not
categorical \textit{because } syntactically incomplete.

As noted above, the full understanding of the Löwenheim-Skolem Theorem
showed, after a few years, that categoricity is impossible for \textit{all}
the formal (or, more generally, with a semantically complete language)
Systems equipped with at least one infinite model (the case of the
ordinary Disciplines). So, this fact is true \textit{regardless if
the Theory is} \textit{syntactically complete or not.} Thus, the subject
matter was finally deciphered by consequences of the Completeness
Theorem (from which, in fact, the same Löwenheim-Skolem Theorem derives).
Neither Gödel, or other alert logician, had ever any good reason to
return to a phrase that, ultimately, had diverted, at least, about
the consequences of the First incompleteness Theorem. Which certainly
were very deep, but concerning essentially different features. The
most significant related to the new and disruptive concept of \textit{machine},
as it was elucidated mainly by Church, Turing and subsequently Chaitin.

The sad postscript is that, still today, are frequent claims that,
in essence, ratified the Gödel's unfortunate sentence without any
correction. For example: \texttt{<\textcompwordmark{}<}\textit{the
syntactic incompleteness of the first order arithmetic causes the
semantics incompleteness of the second order} \textit{logic}\texttt{>\textcompwordmark{}>}%
\footnote{ E. Moriconi, \textit{ibidem} (footnote n. 6), traslated by the author.
A similar phrase is repeated in the \textit{abstract} of F. Berto,
\textit{Gödel's first theorem}, ed. Tilgher Genova, fasc. Epistemologia
27, n.1 (2004). %
}. Overlooking the terminology of \texttt{\textquotedbl{}}expressive
order\texttt{\textquotedbl{}} (that is in effect ambiguous, as we
will try to show later), these words seem, firstly, to suggest the
automatic transmission of syntactic incompleteness to the expanded
System, ie to the second order one (first mistake). Secondly, from
syntactic incompleteness and categoricity it is deduced the semantic
incompleteness (second lapse: it would sufficient the categoricity
plus the infinity of a model).

On the other hand, neither for the language of the integral Theory
of reals (at second order, in its original version), informal and
categorical System, it is valid the semantic completeness Theorem.
Despite that its formal version, expressed at first order, is syntactically
complete (as was shown by Tarski). Evidently, the semantic incompleteness
of the language of this Theory is just due to its categoricity plus
the infinity of the standard model.\textbf{ }Here, what role should
be played by the syntactic incompleteness of \textit{PA}?

Even in a more recent book, it is proposed to deduce semantic incompleteness
of the \texttt{\textquotedbl{}}second-order logic\texttt{\textquotedbl{}}
either exploiting the First incompleteness Theorem or the Church-Turing
Theorem%
\footnote{ C. L. De Florio, \textit{Categoricità e modelli intesi}, ed. Franco
Angeli (2007).%
}. In the first case the theorem is illegally applied to the second-order
Arithmetic. In the second case, the author - agree with others scholars
- founds the proof on the fact that \texttt{\textit{<\textcompwordmark{}<}}\textit{if
}{[}by contradiction{]}\textit{ the second-order logic were }{[}semantically{]}\textit{
complete, then there should be a }{[}effective{]}\textit{ procedure
of enumeration of the valid formulas for the second order...} \texttt{>\textcompwordmark{}>}%
\footnote{ C. L. De Florio, \textit{op. cit}., p. 54, traslated by the author.%
}. However, this consequence assumes that the set of the second order
formulas is enumerable. But now, if with the expression \texttt{\textquotedbl{}}second-order
logic\texttt{\textquotedbl{}} we decide to indicate a countable System,
then the categorical Arithmetic is not included in it!

Indeed, until syntactic incompleteness (or completeness) of second
order Arithmetic remains unproven, we do not distinguish alternatives
to use of the Lowenheim-Skolem Theorem to demonstrate the semantic
incompleteness of its language.

Finally, it is often repeated also the old slip that the non-categoricity
of the formal Arithmetic, ie \textit{PA}, is\textit{ caused} by its
syntactic incompleteness. That is sort of like saying that the infinity
of polygons is \textit{caused} by the infinity of isosceles triangles.

\section*{3. Chaitin's licenses}

In 1970, Gregory Chaitin formulated an interesting informatic version
of the First incompleteness Theorem. In its most simple presentation,
it states that any machine that verify the Church-Turing Thesis%
\footnote{ This famous \texttt{\textquotedbl{}}Thesis\texttt{\textquotedbl{}}
is just the assumption that all the machines are logically reproducible
using \textit{recursive functions} and vice versa. The \textit{recursive
functions} representing all possible arithmetic field of calculability.
More accurate and detailed explanations can be found in any good book
on Logic.%
}can proof the \textit{randomness} of a necessarily finite number of
symbolic strings. The \textit{randomness} of a symbolic string is
defined by the impossibility, for the machine itself, to compress
it beyond a certain degree (which has to be fixed before). For any
machine, exists an infinite and predominant number of random strings:
in fact it is easy to conclude that the probability, for a fixed machine,
to compress an arbitrary finite string of symbols is always quite
low (except to consider a really trivial compression degree). Nonetheless,
compressible strings remain certainly infinite; furthermore has to
be emphasized that any ordinary human creation, still encoded in symbols,
is almost always non-random%
\footnote{ Hence the success of file compression techniques known as \texttt{\textquotedbl{}}loss-less\texttt{\textquotedbl{}},
ie without loss or corruption of data, such as \texttt{\textquotedbl{}}zip\texttt{\textquotedbl{}},
\texttt{\textquotedbl{}}tar\texttt{\textquotedbl{}} and so on. However,
for products with superior randomness, such as videos and music, compression
methods can be effective only with loss or alteration of data: the
formats \texttt{\textquotedbl{}}jpeg\texttt{\textquotedbl{}} and \texttt{\textquotedbl{}}mp3\texttt{\textquotedbl{}}
are examples. %
}.

Thanks to the Chaitin's interpretation we know that any machine can
indicate only a finite number of strings that itself cannot compress.
As a result, no compression program, always stopping, can be certainty
stated as ideal, ie not improvable: by contradiction, it would be
able to determine the randomness of any random string, by virtue of
failing to compress it; in violation of this incompleteness Theorem
version%
\footnote{ G. Raguní, \textit{op.} \textit{cit.}, p. 273.%
}.

Unfortunately Chaitin has released many superficial statements, often
uneven, which produce dangerous confusion about the incompleteness
subject, already in itself not so easy. The fact that these conclusions
are valid also for \textit{universal machines}%
\footnote{ A machine is said \textit{universal }if it is able to reproduce logically
the behavior of any other machine. Its existence derives, by the Church-Turing
Thesis, from the existence of\textit{ universal} recursive functions.
An example of universal machine is any \textit{PC}, however modest,
but with unlimited expandable memory.%
}, has led him, first, to neglect that the definition of \textit{randomness}
always, in any case, is referred with respect to a particular fixed
machine. Moreover, for the obvious fact that, on any ordinary computer,
the natural numbers are represented by symbolic strings, he has precipitously
assigned the \textit{randomness} property to the natural numbers.
As a result of this lightness, Chaitin has repeatedly proclaimed to
have discovered the \texttt{\textquotedbl{}}randomness in Arithmetic\texttt{\textquotedbl{}}%
\footnote{ Two examples: \texttt{\textquotedbl{}}\textit{I have shown that there
is randomness in the branch of pure mathematics known as number theory.
My }work\textit{ indicates that - to borrow Einstein's metaphor -
God sometimes plays dice with whole numbers!}\texttt{\textquotedbl{}}
\textit{Randomness in Arithmetic}, Scientific American 259, n. 1 (July
1988); \texttt{\textquotedbl{}}\textit{In a nutshell, Gödel discovered
incompleteness, Turing discovered uncomputability, and I discovered
randomness}\texttt{\textquotedbl{}}, preface to the book \textit{The
unknowable}, ed. Springer-Verlag, Singapore (1999). This kind of frases,
however, is repeated in almost all his publications, including the
most recent.%
}.

We reiterate that the \textit{randomness} property (originally defined
by A. Kolmorov) only concerns the \emph{strings} of characters and
affects the natural numbers just for the \textit{code} chosen for
them. Encoding is totally arbitrary by the point of view of the Logic.
In effect it is possible, in a specific machine, to use a code that,
although unquestionably uncomfortable and expensive in \textit{bits},
makes \textit{finite} the amount of random natural numbers (or, more
exactly, of the strings which represent them)%
\footnote{ Two examples in G. Raguní, \textit{op. cit.}, p. 276 and 281-282.%
}. Furthermore, as before stated, the same randomness of a symbolic
string is not absolute, but relative with respect to the code and
inner working of the prefixed machine. Given an arbitrary and long
enough string that is random for a particular machine, there is nothing
to ban the existence of another machine, possibly designed \textit{ad
hoc}, which gets to compress it. For that machine, predictably, some
strings that are compressible for the first machines, will be instead
random.

In the same article of \textit{Scientific American} cited in the penultimate
footnote, Chaitin writes:
\begin{quotation}
\textit{How have the incompleteness theorem of Gödel, the halting
problem of Turing and my own work affected mathematics? The fact is
that most mathematicians have shrugged off the results }{[}...{]}\textit{.
They dismiss the fact as not applying directly to their work. Unfortunately
}{[}...{]}\textit{ algorithmic information theory has shown that incompleteness
and randomness are natural and pervasive. This suggests to me that
the possibility of searching for new axioms applying to the whole
numbers should perhaps be taken more seriously.}

\textit{Indeed, the fact that many mathematical problems have remained
unsolved for hundreds and even thousands of years tends to support
my contention. Mathematicians steadfastly assume that the failure
to solve these problems lies strictly within themselves, but could
the fault not lie in the incompleteness of their axioms? For example,
the question of whether there are any perfect odd numbers has defied
an answer since the time of the ancient Greeks. Could it be that the
statement }\texttt{\textit{\textquotedbl{}}}\textit{There are no odd
perfect numbers}\texttt{\textit{\textquotedbl{}}}\textit{ is unprovable?
If it is, perhaps mathematicians had better accept it as an axiom.}

\textit{This may seem like a ridiculous suggestion to most mathematicians,
but to a physicist or a biologist it may not seem so absurd. }{[}...{]}
\textit{Actually in a few cases mathematicians have already taken
unproved but useful conjectures as a basis for their work. }
\end{quotation}
These words seem to suggest the \texttt{\textquotedbl{}}gödelian revolution\texttt{\textquotedbl{}}
of a \texttt{\textquotedbl{}}new\texttt{\textquotedbl{}} Mathematics,
empirical type, which really has always existed: one which makes use
of conjectures. To consider these last as axioms without no metamathematic
justification would be obviously unwise as well as useless. And it
does not sound like a progress but like a resigned presumption: it
desists, a priori, to search for meta-demostrations which often have
been precious sources of development for Logic and Mathematics. \textit{Indeed,
for no undecidable formula, the incompleteness Theorem impedes the
possibility to distinguish it by a purely metamathematic reasoning}.
This erroneous view is repeated, with enthusiasm, in almost all his
recent work: in fact, on the (epidermal, certainly not logic) basis
of the incompleteness Theorem, he goes so far as to question the same
opportunity of the axiomatic Systems%
\footnote{G. Chaitin, \textit{The halting probability $\Omega$: irreducible
complexity in pure mathematics,}

Milan Journal of Mathematics n. 75 (2007), p. 2 y ss. %
}!

The swedish Logician Torkel Franzén, died recently, exposed other
mistakes of Chaitin in 2005%
\footnote{ T. Franzén, \textit{Gödel's Theorem: an incomplete guide to its use
and abuse}, AK Peters (2005).%
}. Here, we refer just one. In the abstract of an article, Chaitin
states: \texttt{<\textcompwordmark{}<}\textit{Gödel's theorem may
be demonstrated using arguments having an information-theoretic flavor.
In such an approach it is possible to argue that if a theorem contains
more information than a given set of axioms, then it is impossible
for the theorem to be derived from the axioms}\texttt{>\textcompwordmark{}>}%
\footnote{ G. Chaitin, \textit{Gödel's Theorem and Information}, International
Journal of Theoretical Physics, n. 22 (1982). %
}. The phrase incorrectly confuses the property that for any machine
there are always infinite random strings, with the proving capability
of the machine itself (which is subject to the incompleteness Theorem).
Franzén confutes the assertion in a simple and irrefutable way: by
the single axiom \texttt{\textquotedbl{}}$\forall$\textit{x}(\textit{x}=\textit{x})\texttt{\textquotedbl{}},
having constant complexity, one can obtain a theorem, type \texttt{\textquotedbl{}}\textit{n}
= \textit{n}\texttt{\textquotedbl{}}, having arbitrarily large complexity
by increasing the number \textit{n.} Indeed, that is guaranteed if
the natural numbers are encoded by any usual exponential code.

The constant \textit{$\Omega$,} introduced by Chaitin with respect
to a fixed universal machine, represents the probability that a random
routine of the machine stops. Its unquestionable interest lies in
the fact that it represents a kind of \textit{best compression} of
mathematical knowledge: knowledge of the first \textit{n} \textit{bits}
of $\Omega$ can solve the halting problem for all programs of length
less than or equal to \textit{n.} However, Chaitin greatly exaggerates
its importance, as far as to describe $\Omega$ as \textit{the way}
to obtain the mathematical knowledge%
\footnote{ For example, in the article: \textit{Meta-mathematics and the foundations
of Mathematics}, EATCS Bulletin, vol. 77, pp. 167-179 (2002), he exposes
a method that in principle\textit{ }should be able to solve the Riemann
hypothesis from knowing a sufficient number of bits of \textit{$\Omega$}.
But the argument is a clear vicious circle.%
}. Naturally, the number $\Omega$ is just a fantastic, insuperable
way to \textit{summarize} this knowledge. After having obtained it
by the traditional theorems and meta-theorems.

This criticism is not meant to attack the figure of this great informatic
logician, but just to clarify a picture that shows still quite confused,
not only to non-experts.

\section*{4. A tired classification}

We now wish to criticize the current classification of the classical
axiomatic Theories in function of its expressive order: \textit{first
order},\textit{ second order}, etc..%
\footnote{A language is said of the first-order if the quantifiers ``$\exists$''
and ``$\forall$'' only may refer to variables (as in the phrase
\texttt{\textquotedbl{}}\textit{every} line that is parallel to line
\textit{r} also is parallel to line \textit{t}\texttt{\textquotedbl{}}).
At second order, also it can be quantified on the predicates (and
to translate phrases like: \texttt{\textquotedbl{}}\textit{Every}
relationship that exists between the lines \textit{r} and \textit{s}
also exists between the lines \textit{r} and \textit{t}\texttt{\textquotedbl{}}).
At third order it can be quantified also on \textit{super-predicates}
(relations between predicates) and so on to infinity.%
}. Unfortunately, it has been consolidated over the years the view
that the formality of language (or, more generally, its semantic completeness)
is a prerogative of the first expressive order and, moreover, that
higher-order language are all necessarily semantic (typically, uncountable).
The error is mainly due to two misunderstandings.

The first is linked to the meaning of \textit{First Order Classical
Logic.} This expression usually refers to the collection of all the
\textit{Classic Predictive Calculi.} Each \textit{Classical Predictive
Calculus} is a first order formal Theory that has got: a) some basic
first order axioms, specified for the first time by Russell and Whitehead,
which formalize the concepts of \texttt{\textquotedbl{}}not\texttt{\textquotedbl{}},
\texttt{\textquotedbl{}}or\texttt{\textquotedbl{}} and \texttt{\textquotedbl{}}exists\texttt{\textquotedbl{}}%
\footnote{ Other usual classical concepts, such as \texttt{\textquotedbl{}}and\texttt{\textquotedbl{}},
\texttt{\textquotedbl{}}imply\texttt{\textquotedbl{}} and \texttt{\textquotedbl{}}every\texttt{\textquotedbl{}},
are defined by them. The first two can be unificated by an unique
connective like NAND (or NOR, Sheffer 1913).%
}; b) some other own axioms, again enumerable and at first order, which
formalize some particular concepts that have to be used in the Theory
(eg, \texttt{\textquotedbl{}}equal\texttt{\textquotedbl{}}, \texttt{\textquotedbl{}}greater
than\texttt{\textquotedbl{}}, \texttt{\textquotedbl{}}orthogonal\texttt{\textquotedbl{}},
etc.); c) the four classical deductive rules: \textit{substitution,
particularization, generalization} and \textit{modus ponens.} Since
these rules consist of purely \textit{syntactic} operations on axioms
and/or theorems (ie are applied, mechanically, just to symbolic code
of the statements), in any Predictive Calculus - and so in the whole
\textit{Classical First Order Logic} - the formality is always verified.
However, this fact does not mean that \textit{every} first order Theory,
founded on a particular Classical Predictive Calculus and deducting
only by the four classic rules, must be formal. Nothing prevents,
for example, that a Theory add an uncountable number of own axioms
to this Calculus, although all expressed at first-order language:
it would result an intrinsically semantic and therefore not formal
language. The\textit{ First Order Classical Logic}, in other words,
does not include all the classical first order Theories. Infinitely
many of them use an informal and/or semantically incomplete language%
\footnote{In agree: M. Rossberg, \textit{First-Order Logic, Second-Order Logic,
and Completeness}, Hendricks et al. (eds.) Logos Verlag Berlin (2004),
on the \textit{Web}:

http://www.st-andrews.ac.uk/\textasciitilde{}mr30/papers/RossbergCompleteness.pdf%
}.

The second mistake is related to the Lindström's Theorem. This states
that any classical Theory expressed in a semantically complete language
\textit{can} be formulated by a first order language. The disorder
stems to confuse \texttt{\textquotedbl{}}can\texttt{\textquotedbl{}}
with \texttt{\textquotedbl{}}must\texttt{\textquotedbl{}}. \textit{The
theorem does not forbid the semantic completeness, or formality, of
higher-order languages}. It only states that when you have just this
case, the Theory can be re-expressed more easily at first order language.
Unquestionably, this property distinguishes the particular importance
of the first order expressions. A property, on the other hand, that
already is evident thanks to the expressive capability of the formal
Set Theory: any \textit{formal} System, in fact, since is fully representable
inside this Theory - which is of first order - is expressible at the
first order. But this importance should not be radicalized.

Surely, the evident fact that the second order languages are typically
uncountable, and therefore intrinsically semantic, aggravates the
situation. That happens because, if the model is infinite, in the
most general case the predicates vary within a certainly not enumerable
set. But nothing preclude that the axioms limit this variability to
a countable subset and that, in particular, the formality is respected%
\footnote{ In agree: HB Enderton, \textit{Second order and Higher order Logic},
Standford Encyclopedia of Philosophy (2007), on the \textit{Web}:
http://plato.standford.edu/entries/logic-higher-order/.%
}. A concrete example is the System obtained from \textit{PA,} expressing
the \textit{partial} induction principle, ie limited to formulas with
at least one free variable, by a single symbolic formula (instead
of by a meta-mathematic axiomatic scheme). It results a second-order
axiom which still has to be interpreted semantically, since the premises
of the Theory do not specify any syntactic deduction with second order
formulas. However, this semantics is not \textit{intrinsic}, ie ineliminable:
if we represent the System inside the formal Set Theory, this axiom
is translated to a set-theoretic axiomatic scheme that generates an
\textit{enumerable} number of formal inductive axioms. Therefore,
in the same original Theory, the formality could be restored by adding
the appropriate premises with which syntactically operate from the
second order induction axiom, so as to produce the required theorems
about the symbolically declarable properties. But, naturally, there
are strategies more simple to reconstitute the formality%
\footnote{ G. Raguní, \textit{op. cit.}, p. 160-161.%
}.

As a result of this confusion, the non-formal nature of the second
(or more) order Theories is often criticized, on the base of intrinsic
semanticity in \textit{some} of these. And others scholars, rather
than highlight that the problem does not lie in the kind of expressive
order, but in the particular nature of the premises of the Theory,
contest that \texttt{\textquotedbl{}}also some second order Systems
have a formal appearance like those of the first order\texttt{\textquotedbl{}}.
In any case, like \textit{some} of the first order!

In conclusion, the cataloging of classical axiomatic Theories in function
of its language order, in general, misleads about their basic logical
properties. Those are consequence of the structure of the premises,
whereas the language order not always plays a decisive role. The main
instrument of classification remains the respect of the Hilbert's
formality or, more generally, of the semantic completeness.

\section*{5. About the interpretation of the Second incompleteness Theorem}

Finally, we have to disrupt the usual interpretation of the Second
incompleteness Theorem. In reference to a Theory that satisfies the
same hypothesis of the First incompleteness Theorem, the Second one
generalizes the undecidability to a class of statements which, interpreted
in the \textit{standard model} mean \texttt{\textquotedbl{}}this System
is consistent.\texttt{\textquotedbl{}} Its complex demonstration,
only outlined by Gödel, was published by Hilbert and Bernays in 1939.

The usual interpretation of this Theorem, object of our criticism,
is that \texttt{\textquotedbl{}}every Theory that satisfies the hypotheses
of the First incompleteness Theorem can not prove its own consistency.\texttt{\textquotedbl{}}
It seems clear, in fact, that the conclusion that a Theory can not
prove its own consistency is valid for \textit{all the classical Systems},
including non-formal! Moreover, this conclusion does not correspond
to the Second incompleteness Theorem, but to a new metatheorem.

Consider an arbitrary classical System. If it is inconsistent, it
is deprived of models and therefore of any reasonable interpretation
of any statement%
\footnote{ More in depth: of any interpretation respecting the principles of
contradiction and excluded third.%
}. Therefore, only the admitting that a given statement of the Theory
means something, implies agreeing consistency. And indubitably this
also applies if the interpretation of the statement is \texttt{\textquotedbl{}}this
System is consistent\texttt{\textquotedbl{}}.

So, if there is no assurance about the consistency of the Theory (which,
to want to dig deep enough, applies to any mathematical Discipline)
we can't be certain on any interpretation of its language. For example,
in the case of the usual Geometry, when we prove the pythagorean Theorem,
what we really conclude is \texttt{\textquotedbl{}}if the System supports
the Euclidean model (and therefore is consistent), then in every rectangle
triangle c\textsubscript{1}\textsuperscript{2}+c\textsuperscript{2}\textsubscript{2}=I\textsuperscript{2}\texttt{\textquotedbl{}}.
Certainly, a deduction with undeniable epistemological worth, still
in the catastrophic possibility of inconsistence.

But now let's see what's happen if a certain theorem of a certain
Theory is interpreted with the meaning: \texttt{\textquotedbl{}}this
System is consistent\texttt{\textquotedbl{}} in a given interpretation
\textit{M.} Similarly, what we can conclude by this theorem is really:
\texttt{\textquotedbl{}}if the System supports the model \textit{M}
(and therefore is consistent), then the System is consistent.\texttt{\textquotedbl{}}
\textit{Something that we already knew and, above all, that does not
demonstrate at all the consistency of the System.} Unlike any other
statement with a different meaning in \textit{M}, for this kind of
statement we have a peculiar situation:\textit{ bothering to prove
it within the Theory is epistemologically irrelevant in the ambit
of the interpretation M}. In more simple words, the statement in question
can be a theorem or be undecidable with no difference for the epistemological
view. Just it cannot be the denial of a theorem, if \textit{M} is
really a model. So, in any case, \textit{the problem of deducing the
consistency of the Theory is beyond the reach of the Theory itself.}
We propose to call \textit{Metatheorem of undemonstrability internal
of consistency} this totally general metamathematic conclusion.

Then, the fact that in a particular hypothetically consistent Theory,
such a statement is a theorem or is undecidable, is depending on the
System and on the specific form of the statement. For Theories that
satisfy the assumptions of the First, Second incompleteness Theorem
guarantees that \texttt{\textquotedbl{}}normally\texttt{\textquotedbl{}}
these statements are undecidable. We say \texttt{\textquotedbl{}}normally\texttt{\textquotedbl{}}
because apparently there are also statements, still expressing consistency
of the System in other peculiar models that, on the contrary, turn
out to be \textit{theorems} of the Theory%
\footnote{ See eg: G. Lolli, \textit{Da Euclide a Gödel}, ed. Il Mulino (2004),
p. 140 y 142 and A. Martini, \textit{Notazioni ordinali e progressioni
transfinite di teorie}, Thesis Degree, University of Pisa (2006),
p. 11-15, on the \textit{Web}: http://etd.adm.unipi.it/theses/available/etd-11082006-161824/unrestricted/tesi.pdf%
}. As the same Lolli says, \texttt{\textquotedbl{}}it seems that not
even a proof shuts discussions\texttt{\textquotedbl{}}%
\footnote{ G. Lolli, \textit{op. cit.}, p. 142.%
}. But in any case, as before concluded, this debate cannot affect
the validity of the proposed Metatheorem of undemonstrability internal
of consistency.

In conclusion, the Second incompleteness Theorem identifies another
class of essentially undecidable statements for any Theory that satisfies
the hypothesis of the First one. Whilst the First incompleteness Theorem
determines only the Gödel's statement, the Second extends the undecidability
to a much broader category of propositions. But, contrary to what
is commonly believed, this drastic generalization does not introduce
any new and dramatic epistemological concept about the consistency
of the System. It doesn't, even if the Theorem were valid \textit{for
every} statement interpretable as \texttt{\textquotedbl{}}this System
is consistent\texttt{\textquotedbl{}} (which, we reaffirm, seems to
be false). Because by it, in any case, we cannot conclude that \texttt{\textquotedbl{}}the
System cannot prove its own consistency\texttt{\textquotedbl{}}: this
judgment belongs to a completely general metatheorem which seems never
have been stated, despite its obviousness and undeniability%
\footnote{ The consistency of a System may be demonstrated only outside the
same, by another external System. For which, in turn, the problem
of consistency arises again. To get out of this endless chain, the
\texttt{\textquotedbl{}}last\texttt{\textquotedbl{}} conclusion of
consistency has to be purely metamathematic. Actually, the most general
Theory (that demonstrates the consistency of the ordinary mathematic
Disciplines) is the formal Set Theory and the conclusion of its consistency
only consists in a \texttt{\textquotedbl{}}reasonable conviction\texttt{\textquotedbl{}}.%
}.

Often the misreading is aggravated for a sort of incorrect \texttt{\textquotedbl{}}intuitive
proof\texttt{\textquotedbl{}} of the Second incompleteness Theorem;
that sounds like: \texttt{\textquotedbl{}}Let \textit{S }be a System
that satisfies the hypotheses of the First incompleteness Theorem
and \textit{C} its statement affirming the consistency of \textit{S.}
The First incompleteness Theorem shows that if \textit{S} is consistent,
the Gödel's statement, \textit{G}, is undecidable. Now, if \textit{C}
were provable, we could deduce that \textit{G} is undecidable and
therefore unprovable. But since \textit{G} claims to be itself unprovable,
this fact would mean just to prove \textit{G}, that is absurd. Therefore,
\textit{C} must be unprovable\texttt{\textquotedbl{}}%
\footnote{ See eg: P. Odifreddi,\textit{ Metamorfosi di un Teorema}, (1994),
on the \textit{Web}: http://www.vialattea.net/odifreddi/godel.htm%
}. The flaw is obvious: the reasoning gives to \textit{C} and \textit{G}
a semantic value which is certified only assuming that the System
admits a model with such interpretations and therefore already assuming
that it is consistent. \textit{In this model} there is no doubt that
the truth of \textit{C} implies the truth of \textit{G}, but the syntactic
implication \textit{C→G }is a totally different question. In general,
the possibility that \textit{C} is a theorem causes no absurd consequence,
because that does not imply really the consistency of the System,
as it was observed. On the other hand, in case of inconsistency, maybe
doesn't happen that every statement is a theorem?

Actually, the syntactic inference \textit{C→G} is not so trivial and,
as noted, it does not apply in all cases but depends on the symbolic
structure of the statement \textit{C}.

The only epistemological result of effective importance about the
consistency is due to the Metatheorem of undemonstrability internal
of consistency. And we emphasize that its meta-demonstration, since
refers to any arbitrary classical System, must consist in a purely
meta-mathematical reasoning (like that one we did), ie it cannot be
formalized.

\begin{center}
-------- 
\par\end{center}

In the book, already mentioned, \textit{Confines Lógicos de la Matemática}%
\footnote{ Also in italian: \textit{I confini logici della matematica}, on ed.
Aracne, Bubok, Scribd, Lulu or Amazon. %
}, we propose other revisions and some proposals to update the terminology
of these arguments, which continues unchanged since the time of Hilbert.
\end{document}